\documentclass{article}
\usepackage{german}
\usepackage[latin1]{inputenc}
\usepackage{enumerate}
\usepackage{longtable}
\usepackage{graphicx}
\usepackage{fnpara}
\usepackage{amsmath}
\usepackage{makeidx}
\usepackage{a4wide}
\usepackage{color}
\usepackage{amsopn}
\usepackage{ulem}
\usepackage{float}
\usepackage{yfonts}
\graphicspath{{EPSBilder/}{PDFBilder/}}
\newtheorem{defi}{Definition}[section]

\newtheorem{anm}{Anmerkung}[section]

\newtheorem{satz}{Satz}[section]
\begin{document}
\title{Zur iterativen Lösung von linearen Gleichungssystemen}
\author{H. Karl und S. Karl,  GSO - Hochschule Nürnberg}
\date{1. 12. 2012}
\maketitle

\begin{abstract}
Notwendiges und hinreichendes Kriterium dafür, dass eine Fixpunktiteration vom Typ
\begin{align*}
\underline{x}_{m+1}= \underline{\Phi}\ \cdot  \underline{x}_{m}+\underline{h};
\end{align*}

\begin{center}   
 \( \qquad \qquad \qquad \qquad \qquad \qquad \qquad \qquad \qquad \qquad \qquad
\underline{\Phi} \in \mathbf{R}^{n \times n};\   \underline{x}_m, \underline{h} \in \mathbf{R}^{n};\  m \in \mathbf{N}_{0},n \in \mathbf{N}
\)
\end{center}
unabhängig vom Anfangswert gegen einen Fixpunkt \( \underline{x}^* \in \mathbf{R}^{n}\) konvergiert  (für  den  dann: 
 \(\underline{x}^*=\underline{\Phi} \  \underline{x}^* + \underline{h} \) erfüllt ist) ist  bekanntlich, dass für alle Eigenwerte von \( \underline{\Phi}: \ \lambda_1,..., \lambda_l \ ( l \le n) \) gilt: \( \mid \lambda_{i} \mid < 1. \text{ \quad  i} \in \lbrace 1,...,l \rbrace  \)
\newline \newline In dieser Arbeit wird ein neues Verfahren angegeben, das es gestattet die Konvergenz der Fixpunktiteration auch für den Fall zu erzwingen, dass die Eigenwerte von \(\underline{\Phi}\) (der ''Iterationsmatrix'')  betragsmäßig größer als eins sind.
\newline
\newline Im ersten Kapitel dieser Arbeit wird ein neues Konzept vorgestellt, das es erlaubt, eine vorgelegte Fixpunktiteration so zu erweitern, dass jedes Eigenwerttupel der neu entstehenden Iterationsmatrix eingestellt werden kann. Das Eigenwerttupel der ursprünglichen Iterationsmatrix \(\underline{\Phi}\) muss aber noch bekannt sein. Ein notwendiges und hinreichendes Kriterium für die Durchführbarkeit dieses Konzepts wird angegeben. Aus den errechneten Ergebnissen kann sehr einfach auf die Lösung des Ausgangsproblems zurückrechnet werden.
\newline 
\newline Im zweiten Kapitel wird dieses Konzept dann so erweitert, dass man auf die Kenntnis des Eigenwerttupels von \(\underline{\Phi}\) verzichten kann: Die Verwendung von Optimierungsverfahren für lineare dynamische Systeme erzwingt eine Iterationsvorschrift, bei welcher der spektrale Radius der Iterationsmatrix kleiner eins (aber ansonsten nicht genauer bekannt) ist.
\newline
\newline Das dritte Kapitel liefert schließlich eine Iterationsvorschrift, bei der der spektrale Radius vorgeschrieben werden kann, ohne dass das Spektrum der Iterationmatrix \(\underline{\Phi}\) bekannt sein muß. Insgesamt wird sich zeigen, dass sich damit das Einsatzgebiet von iterativen Rechenverfahren zur Lösung von linearen Gleichungssystemen bedeutend erweitern läßt.
\end{abstract}

\newpage
\section{Ein neues Konzept für die iterative Berechnung der Lösung eines linearen Gleichungssystems}
Für die Überführung eines linearen Gleichungssystems vom Typ 
\[ \qquad \qquad \underline{A} \cdot \underline{x}  = \underline{b};\qquad  \underline{A} \in \mathbf{R}^{n \times n}, \quad \underline{x}, \underline{b} \in \mathbf{R}^{n} \]
in die Matrixdifferenzengleichung
\begin{align*} 
\qquad \quad \underline{x}_{m+1}= \underline{\Phi}\ \cdot \underline{x}_{m}+\underline{h}; \qquad m \in \mathbf{N}_{0},\ \underline{x}_{0}=\mathbf{R}	^{n} 
\end{align*}
wird hier beispielhaft das Jacobische Iterationsverfahren verwendet. Aber auch für andere Iterationsverfahren (wie z.B. das Gauß-Seidelsche Verfahren) bleibt das hier vorzustellende Lösungskonzept gültig.

\subsection{Eine Erweiterung des Jacobischen Iterationsverfahrens}
 Ausgehend von 
\begin{align}
\underline{A} \cdot \underline{x} = \underline{b}\ \ ;  
\qquad \qquad \underline{A} \in \mathbf{R}^{n \times n} \qquad \underline{b} , \underline{x} \in \mathbf{R}^n
\end{align}
erhält man bekanntlich nach Einführung von:

 \(   \qquad \qquad \underline{L} = (l_{ij})_{i,j = 1,...,n} 
\qquad \  l_{i,j}= \left\{ \begin{array}{r@{\quad:\quad}l}
                       a_{i,j} & i>j \\
                       0 & \text{sonst}
                        \end{array} \right.
  \)
  
  \(   \qquad \qquad \underline{U} = (u_{ij})_{i,j = 1,...,n} 
\qquad u_{i,j}= \left\{ \begin{array}{r@{\quad:\quad}l}
                       a_{i,j} & i<j \\
                       0 & \text{sonst}
                        \end{array} \right.
  \)
  
  \(   \qquad \qquad \underline{D} = (l_{ij})_{i,j = 1,...,n} 
\qquad \  d_{i,j}= \left\{ \begin{array}{r@{\quad:\quad}l}
                       a_{i,j} & i=j \\
                       0 & \text{sonst}
                        \end{array} \right.
  \)
\newline für (1) die Darstellung:
\begin{align}
(\underline{L}+\underline{D}+\underline{U})\underline{x} =\underline{b}
\end{align}
und weiter (wenn \(\underline{D}\) invertierbar ist)

\begin{align} \underline{x}=-\underline{D}^{-1} \cdot(\underline{L}+\underline{U})\cdot\underline{x}+\underline{D}^{-1}\cdot\underline{b}
\end{align}
Kürzt man ab: 
\begin{align}
\underline{\Phi} &=-\underline{D}^{-1}\cdot(\underline{L}+\underline{U}) \\
\underline{h}&=\underline{D}^{-1}\cdot\underline{b} 
\end{align} 
so kommt man auf die zugehörige Fixpunktform
\begin{align}
\underline{x}= \underline{\Phi}\ \cdot  \underline{x}+\underline{h} 
\end{align}
Davon ausgehend wird die Iterationsvorschrift (Matrixdifferenzengleichung) 
\begin{align} 
\underline{x}_{m+1}= \underline{\Phi}\ \cdot \underline{x}_{m}+\underline{h};\qquad m \in \mathbf{N}_{0},\underline{x}_{0}=\mathbf{R}	^{n}  
\end{align}
formuliert.
Man definiert:
\begin{defi}
Das Iterationsverfahren (7) ist konsistent zur Matrix \underline{A}, wenn für alle \underline{b} \( \in \mathbf{R}^{n} \) der Term \( \underline{A}^{-1}\cdot \underline{b} \) ein Fixpunkt für dieses Iterationsverfahren ist.
\end{defi}
Bekanntlich ist die Iterationsvorschrift (7) konsistent zur Matrix \underline{A}, wenn für den Spektralradius von \( \underline{\Phi} \) gilt [Mei] :
\begin{align}
  \rho \left( \underline{\Phi} \right) < 1  
\end{align}
Für den Fall \(   \rho \left( \Phi \right) \ge 1 \) werde die Iterationsvorschrift (7) jetzt wie folgt erweitert:
\begin{align} 
\underline{\tilde{x}}_{m+1}= 
\underbrace{ \left(\underline{\Phi}-\underline{h}\cdot \underline{k}^{T} \right)}_{=\underline{\tilde{\Phi}}}
\cdot\underline{\tilde{x}}_{m}+\underline{h};\qquad \underline{\tilde{x}}_{0}=\underline{x}_{0}; \qquad \underline{k} \in \mathbf{R}^{n} 
\end{align}

\begin{anm}
Die Differenzengleichung (9) wird in dieser Arbeit das erweiterte Iterationsverfahren genannt.
\end{anm}

Durch die Forderung nach einem Eigenwerttupel für \( \underline{\tilde{\Phi}} \) soll der Vektor  \( \underline{k}\) festgelegt werden. Diese Forderung beschreibt das (in der Regelungstechnik seit langem verwendete) Konzept der Zustandsregelung. Unter welchen Bedingungen diese Forderung erfüllt werden kann, besagt der nachfolgende Satz [Hart1]:
\begin{satz}
Notwendig und hinreichend dafür, dass für eine Matrix \( \underline{\tilde{\Phi}} = \underline{\Phi}-\underline{h}\cdot \underline{k}^{T} \) durch geeignete Wahl von \underline{k} jedes beliebige Eigenwerttupel eingestellt werden kann ist, dass: \newline 
\begin{center}
\(rg\left(\underline{h},\underline{\Phi} \cdot \underline{h}, \cdots \underline{\Phi}^{n-1} \cdot \underline{h}\right) = n \)
\end{center}
erfüllt ist.
\end{satz}
Rechenverfahren zur Berechnung eines geeigneten \(\underline{k}\) -  Vektors nachdem 
\[ \text{die Eigenwerte } \lambda_{1},\cdots,\lambda_{n}  \text{ der Matrix  } \underline{\Phi} \text{  bekannt} \]und 
\[ \text { die Eigenwerte } \tilde{\lambda}_{1}, \cdots , \tilde{\lambda}_{n} \text{ der Matrix } \tilde {\underline{\Phi}} \text{ gefordert} \] werden, sind seit langem bekannt. Sie laufen auf einen Koeffizientenvergleich der zughörigen charakteristischen Polynome (von \(\underline{\Phi} \) und \(\tilde {\underline{\Phi}} \) ) und einen Basiswechsel im \(\mathbf{R}^{n}\) hinaus. Um aber diejenige Transformationsmatrix zu berechnen, die diesen Basiswechsel vermittelt, ist die Lösung eines weiteren \(n \times n\) - Gleichungssystems erforderlich. [Foell]
\[ \]
Damit erscheint das gesamte, neue Konzept für die iterative Berechnung der Lösung eines linearen Gleichungssystems hinfällig. Der Aufwand ist normalerweise zu groß: damit das Gleichungssystem (1) in der Darstellung (9) iterativ gelöst werden kann, müßte ja vorher ein anderes Gleichungssystem von derselben Ordnung gelöst werden. Hinzu kommt die Berechnung des Eigenwerttupels von \(\underline{\Phi}\). Immerhin sollen einige bemerkenswerte Aspekte für dieses neue Konzept vorgestellt werden. Damit rechtfertigen sich weitere Modifikationen die schließlich auf ein leistungsfähiges Lösungskonzept (vgl. Folgekapitel) führen.
\[ \]
Von den neuen Aspekten seien genannt (diese gelten, wenn die Bedingung von Satz 1.1 erfüllt ist):
\begin{enumerate}
\item Unabhängig vom Spektralradius der Matrix \(\underline{\Phi}\) kann man eine gegen den Fixpunkt von (9) konvergente Iteration erzeugen.
\item Indem man von \(\underline{\tilde{\Phi}}\) verlangt: \( \tilde{\lambda}_{1}=0, \cdots , \tilde{\lambda}_{n}=0\), konvergiert die Iteration (9) nach genau ''n'' - Iterationsschritten und das Ergebnis ist exakt.
\item Die Rückrechnung von dem unter (9) errechneten Fixpunkt in den Fixpunkt von (7) und damit in die Lösung von (1) ist elementar.
\end{enumerate}
Bevor diese neuen Aspekte anhand eines kleinen Rechenbeispiels demonstriert werden, soll der unter Punkt 2 genannte Aspekt hervorgehoben werden.
\subsection{Minimale Anzahl von Iterationen}
Die Bedingung von Satz 1.1 sei erfüllt. \underline{k} kann dann für jedes geforderte Eigenwerttupel von \(\underline{\tilde{\Phi}}\) (also von \(\tilde{\lambda}_{1}, \cdots , \tilde{\lambda}_{n}  \)) berechnet werden. Insbesondere auch für \( \tilde{\lambda}_{1}=0, \cdots , \tilde{\lambda}_{n}=0  \). \(\underline{\tilde{\Phi}}\) ist dann nilpotent. Dieser Fall soll jetzt für die \(\underline{\tilde{\Phi}}\) - Matrix der Differenzengleichung (9) vorliegen: 
\begin{align*} 
\underline{\tilde{x}}_{m+1}= \underline{\tilde{\Phi}} \cdot 
\underline{\tilde{x}}_{m}+\underline{h};\qquad \underline{\tilde{x}}_{0}=\underline{x}_{0} 
\end{align*}
Der zugehörige Fixpunkt sei \(\underline{\hat{x}}\):
\begin{align*} 
\underline{\hat{x}}= \underline{\tilde{\Phi}} \cdot 
\underline{\hat{x}}+\underline{h} \qquad \qquad \quad \quad
\end{align*}
Subtrahiert man die letzten beiden Gleichungen voneinander, so erhält man
\begin{align}
\underline{\tilde{x}}_{m+1}-\underline{\hat{x}}= \underline{\tilde{\Phi}} \cdot 
\left( \underline{\tilde{x}}_{m} -\underline{\hat{x}} \right);\qquad \underline{\tilde{x}}_{0}-\underline{\hat{x}}=\underline{x}_{0} - \underline{\hat{x}}
\end{align}
und mit 
\begin{align}
\underline{e}_{m}:= \underline{\tilde{x}}_{m} - \underline{\hat{x}}
\end{align}
endlich
\begin{align}
\underline{e}_{m+1}= \underline{\tilde{\Phi}} \cdot \underline{e}_{m} \ ; \quad \underline{e}_{0}= \underline{x}_{0} - \underline{\hat{x}}
\end{align}
Rückwärtseinsetzen der einzelnen Gleichungen aus (12) liefert
\begin{align*}
\underline{e}_{1}&= \underline{\tilde{\Phi}} \cdot \underline{e}_{0} \\
\underline{e}_{2}&= \underline{\tilde{\Phi}}^{2} \cdot \underline{e}_{0} \\
\vdots \\
\underline{e}_{n}&= \underline{\tilde{\Phi}}^{n} \cdot \underline{e}_{0} 
\end{align*}
Wegen der vorausgesetzten Nilpotenz von \(\underline{\tilde{\Phi}}\) ist \(\underline{e}_{n}=\underline{0}\) und folglich \(\underline{\tilde{x}}_{n}=\underline{\hat{x}} \). Nach \(m = n\) Schritten wird der Fixpunkt also exakt erreicht. 
\newline
\newline
Im anschließenden Beispiel wird noch gezeigt, wie sich aus diesem Fixpunkt  \(\underline{\hat{x}} \) für das erweiterte Iterationsverfahren die Lösung  \(\underline{x} \) für das Ausgangsproblem (1) berechnet.

\subsection{Beispiele}
\subsubsection{Sämtliche Eigenwerte von \(\underline{\tilde{\Phi}}\) liegen bei Null}
Für das entsprechend (1) zu lösende, lineare Gleichungssystem lauten Matrix und Vektor
\begin{align*}
\underline{A} = 
\begin{pmatrix}
59 & -63 & -2 \\
29 & 42 & 51 \\
36 & 31 & -67
\end{pmatrix}
\quad \text{  und } \qquad \underline{b}=
\begin{pmatrix}
-73 \\
266\\
-103
\end{pmatrix}
\quad  \text{,   mit Lösungsvektor    } \qquad \underline{x}=
\begin{pmatrix}
1 \\
2\\
3
\end{pmatrix}.
\end{align*}
Die zugehörigen Matrizen entsprechend (7) sind
\begin{align*}
\underline{\Phi}= 
\begin{pmatrix}
0,000.000.000 & 1,067.796.610 & 0,033.898.305 \\
-0,690.476.190 & 0,000.000.000 & -1,214.285.714 \\
0,537.313.433 & 0,462.686.567 & 0,000.000.000
\end{pmatrix} 
\text{  und  }
\quad \underline{h}=
\begin{pmatrix}
&-1,237.288.135.593 \\
&\ 6,333.333.333.333\\
&\ 1,537.313.432.836
\end{pmatrix}
\end{align*}
Da \(\lambda_1= 0,235.428 +i\cdot 1,202.990, \lambda_2=0,235.428 -i\cdot 1,202.990, \lambda_3 = -0,470.856 \) die Eigenwerte von \(\underline{\Phi}\) (und damit \(|\lambda_{1}|, |\lambda_{2}|>1\)) sind, wird das Iterationsverfahren (7) nicht konvergieren.
\newline
\[\]
\newline
Verlangt man nunmehr, dass \(\underline{\tilde{\Phi}}\) (vgl. (9)) einen Dreifacheigenwert bei \(\quad \tilde{\lambda}_{1}=0,\quad  \tilde{\lambda}_{2}=0,\quad \tilde{\lambda}_{3}=0\) habe, dann ergibt sich (nach Koeffizientenvergleich der charakteristischen Polynome von \(\underline{\Phi}\) und \(\underline{\tilde{\Phi}}\) und anschließendem Basiswechsel im \(\mathbf{R}^{n}\) [Foell]) ein Vektor 
\begin{align*}
\underline{k}^{T} = 
\begin{pmatrix}
-0,117.899 & 0,025.369 & -0,199.404
\end{pmatrix}
\end{align*}
Nach drei Schritten erreicht die Iteration (9) den Vektor

\begin{align*}
\underline{\tilde{x}}_{3} = 
\begin{pmatrix}
 \ \ 2,988.423 \\
\ \ 5,976.846 \\
\  \ 8,956.270
\end{pmatrix}
\end{align*}
(Der Startvektor war in diesem Fall der Nullvektor.) Zufolge 1.2 ist (unabhängig vom Startvektor) \(\underline{\tilde{x}}_{3}\) der Fixpunkt \(\underline{\hat{x}}\) der Iteration (9).
\newline Den originalen Lösungsvektor ''\(\underline{x}\)'' erhält man, wenn man,  ausgehend von (9), den Wert für \(\underline{\hat{x}}\) einsetzt
\begin{align*} 
\underline{\hat{x}}= 
 \left(\underline{\Phi}-\underline{h}\cdot \underline{k}^{T} \right)\cdot
\underline{\hat{x}}+\underline{h}
\end{align*}
auflöst (vgl. (4) und (5))
\begin{align*}
\underline{\hat{x}}&=\Phi \cdot \underline{\hat{x}}-\underline{h}\cdot \underline{k}^{T}\cdot \underline{\hat{x}} + \underline{h} \\
&=-\underline{D}^{-1}\cdot(\underline{L}+\underline{U})\cdot \underline{\hat{x}} -\underline{D}^{-1}\cdot\underline{b}\cdot \underline{k}^{T}\cdot \underline{\hat{x}}+\underline{D}^{-1}\cdot\underline{b}
\end{align*}
mit \(\underline{D}\) auf beiden Seiten der Gleichung durchmultipliziert
\begin{align*}
\underline{D}\cdot \underline{\hat{x}}&=-(\underline{L}+\underline{U})\cdot \underline{\hat{x}}-\underline{b}\cdot \underline{k}^{T}\cdot \underline{\hat{x}}
+ \underline{b}
\end{align*}
umstellt und rücksubstituiert:
\begin{align*}
\underline{A}\cdot \underline{\hat{x}}&=\underline{b}-\underline{b}\cdot \underline{k}^{T}\cdot \underline{\hat{x}}
\end{align*}
Dann erhält man
\begin{align*}
\underline{\hat{x}}&=\underline{A}^{-1}\cdot\underline{b}-\underline{A}^{-1}\cdot\underline{b}\cdot \underline{k}^{T}\cdot \underline{\hat{x}} \\
&= \underline{x}-\underline{x} \cdot \left(\underline{k}^{T} \cdot \underline{\hat{x}}\right)
\end{align*}
und endlich
\begin{align}\
\underline{x}= \left(\frac{1}{1-\underline{k}^{T} \cdot \underline{\hat{x}}}\right) \cdot \underline{\hat{x}}
\end{align}
%
%
Berechnet man schließlich mit diesem Lösungsvektor den Restfehler \(\underline{r}:=\underline{A} \cdot\underline{x} - \underline{b} \), so erhält man für das vorgelegte Beispiel Restfehlerkomponenten von der Größenordnung \(10^{-14}\).
\subsubsection{Die Eigenwerte von \(\underline{\tilde{\Phi}}\) liegen bei \(\tilde{\lambda}_1 =0,1\) , \(\tilde{\lambda}_{2}= 0,2\) und \(\tilde{\lambda}_{3}= 0,3\)}
Hierfür lautet 
\begin{align*}
\underline{k}^{T} = 
\begin{pmatrix}
-0,119.694 & -0,072.227 & -0.189.066
\end{pmatrix}
\end{align*}
Um die Genauigkeit von 1.3.1 zu erreichen sind 45 Iterationsschritte notwendig. Es ist
\begin{align*}
\underline{\tilde{x}}_{45} = 
\begin{pmatrix}
 \ \ 5,929.411 \\
\ \ 11,858.823 \\
\  \ 17,788.235
\end{pmatrix}
\end{align*}
Der Startvektor für dieses Iterationsergebnis ist \(x_{0}=\left(1\ 1\  1\right)^{T}\). Nach Rücktransformation (13) ergibt sich der Lösungsvektor mit Restfehlerkomponenten wieder von der Größenordnung \(10^{-14}\).

\section{Anwendung von Optimierungsverfahren}
Bestimmte Optimierungsverfahren für lineare Regelsysteme liefern ebenfalls einen Vektor ''\(\underline{k}\)'' der einen Spektralradius \(\rho \left( \tilde{\underline{\Phi}} \right) < 1\) erzwingt. Es ist naheliegend, diesen Vektor ''\(\underline{k}\)'' auch für das erweiterte Iterationsverfahren (9) einzusetzen. Diese Optimierungsverfahren benötigen keine Berechnung von Matrixeigenwerten. Die Grundlagen zu dem hier vorzustellenden Optimierungskonzept sind z.B. in [Ludyk] dargestellt. Eine typische Aufgabenstellung und die zugehörige Lösung für die Optimierung linearer, dynamischer Systeme wird jetzt vorgestellt. Anschließend wird der, für das Optimierungsverfahren gefundene, Vektor ''\(\underline{k}\)'' in unserem Iterationsverfahren verwendet.    
\subsection{Das Optimierungsproblem}
Gegeben sei für \(\quad \underline{\Phi} \in \mathbf{R}^{n \times n},\  \underline{h} \in \mathbf{R}^{n}, \lbrace u_{m}\rbrace_{m \in \mathbf{N_{0}}} \  \text{ mit } \ u_{m}\in \mathbf{R} 
\)
\newline
\newline die Iterationsvorschrift für die Folge \( \quad \lbrace \underline{x}_{m}\rbrace_{m \in \mathbf{N_{0}}} \ \text{ mit } \  \underline{x}_{m}\in \mathbf{R}^{n} \):

\begin{align}
\underline{x}_{m+1}= \underline{\Phi}\ \cdot \underline{x}_{m}+\underline{h}\cdot u_{m}\ 
\end{align}

Weiterhin seien gegeben:
\begin{itemize}
\item die positiv definite Matrix \(\underline{S} \in \mathbf{R}^{n \times n}\)
\item die Zahl \(r > 0 \)
\item ein Vektor \( \underline{c} \in \mathbf{R}^{n} \) mit der Eigenschaft, dass gilt:
\end{itemize}

\begin{align}
 \text{Rang}
\left( \begin{array}{c}
\underline{c}^{T} \\
\underline{c}^{T} \cdot \underline{\Phi}\\
\vdots \\
\underline{c}^{T} \cdot \underline{\Phi}^{n-1}
\end{array} \right) = n \qquad
\text{ und } 
\qquad \underline{c}\cdot \underline{c}^{T} = \underline{Q} \in \mathbf{R}^{n \times n} \text{ ist positiv semidefinit}
\end{align}

Gesucht ist für beliebig vorgegebens \( N \in \mathbf{N}\) eine Vektorfolge \(\lbrace \underline{k}_{m}\rbrace _{\ m \in \lbrace 0,1,2,...,N-1\rbrace }\), sodass nach Bildung von

\begin{align} 
\underline{z}_{m+1}&= \underline{\Phi}\ \cdot \underline{z}_{m}+\underline{h}\cdot u_{m};\qquad \underline{z}_{m=0}=\underline{z}_{0}  \\
u_{m}&=-\underline{k}^{T}_{m}\cdot \underline{z}_{m}; 
\end{align}

das quadratische Gütekriterium 
\begin{align}
J_{N} = \frac{1}{2}\ \underline{z}_{N}^{T}\cdot \underline{S}\cdot \underline{z}_{N}+\frac{1}{2} \sum_{m=0}^{N-1} \underline{z}_{m}^{T}\left(\underline{c}\cdot \underline{c}^{T}\right) \underline{z}_{m} + r\cdot u_{m}^{2}\qquad  N\in \mathbf{N}
\end{align}
minimal wird.
\[\]
Der Wechsel der Variablenbezeichnung von ''\(\underline{x}_{m}\)'' nach ''\(\underline{z}_{m}\)'' wird notwendig, da nach Einsetzen von (17) in (16) im Gegensatz zu (9) nunmehr eine homogene Matrix - Differenzengleichung entsteht. Man beachte weiterhin, dass in (16) mit variablem ''\(\underline{k}_{m}\)'' gerechnet wird.
\[\]
Die Lösung dieser Optimierungsaufgabe lautet unter der Voraussetzung, dass zusätzlich Satz 1.1 gilt [Ludyk]:
\begin{satz}
Für jedes N \(\in \mathbf{N} \) gilt: Für das mathematische Modell (16)-(18) in Verbindung mit der Forderung: \[J_{N} = Min \] 
berechnet sich 
der ''optimale Vektor'' zu
\begin{eqnarray}
 \underline{k}_{m}^{T}=\frac{1}{r+\underline{h}^{T}\cdot \underline{P}_{m+1}\cdot \underline{h}}\cdot \left(\underline{h}^{T}\cdot \underline{P}_{m+1} \cdot \underline{\Phi} \right) ; \qquad m\in \lbrace 0,1,2,...,N-1 \rbrace
 \end{eqnarray}
Hierbei ist \( \underline{P}_{m} \in \mathbf{R}^{n \times n} \) die Lösung der Matrix - Riccati - Differenzengleichung
\begin{eqnarray} \underline{P}_{m}=\underline{Q}+\underline{\Phi}^{T}\left[\underline{P}_{m+1}-\underline{P}_{m+1}\underline{h}\underline{h}^{T}\underline{P}_{m+1}\cdot \frac{1}{r+\underline{h}^{T}\underline{P}_{m+1}\underline{h}}  \right]\underline{\Phi}\  ;\quad m\in \lbrace 1,2,...,N-1 \rbrace 
\end{eqnarray}
mit dem Endwert \(\underline{P}_{N}=\underline{S}\).
\end{satz}
Die Matrix - Riccati - Differenzengleichung ist also ''rückwärts'' zu lösen. Aus der so entstehenden Folge \(\lbrace{\underline{P}_{N},\underline{P}_{N-1}, \cdots ,\underline{P}_{1} \rbrace} \) errechnet sich \(\lbrace{ \underline{k}_{N-1},\underline{k}_{N-2},\cdots,\underline {k}_{0} \rbrace}\).
\newline
\newline Der sich für ''\(J_{N}\)'' ergebende Wert ist für unsere Betrachtungen nicht weiter interessant.\newline  Wichtig hingegen ist, daß die damit verbundene Rekursion (die homogene Matrix - Differenzengleichung)
\[\underline{z}_{m+1}= \left(\underline{\Phi}-\underline{h} \cdot \underline{k}_{m}^{T}\right)\underline{z}_{m}; \qquad \underline{z}_{m=0}=\underline{z}_{0} \]
unabhängig vom Anfangsvektor \(\underline{z}_{0}\)
gegen den Grenzwert \(\lim_{N \rightarrow \infty} \underline{z}_{N}=\underline{0}\) konvergiert. (Das ist zumindest plausibel, da wegen der vorausgesetzten positiven Definitheit von \(\underline{S}\) und der Semidefinitheit von \(\underline{Q}\) ein Minimum für ''\(J_{N}\)'' nur angenommen werden kann, wenn \(\lim_{N \rightarrow \infty} \underline{z}_{N}=\underline{0}\) gilt. Man beachte, daß für ''N'' keine weiteren Annahmen als \(N \in \mathbf{N}\) gemacht wurde und dass folglich ''N'' beliebig groß sein darf. Tatsächlich kann man zeigen, dass aufgrund der vorausgesetzten Eigenschaften von Vektor \(\underline{c}\),  vgl. (15), die positive Semidefinitheit von \(\underline{Q}\) hierfür ausreicht.[Ludyk]) 
\begin{anm}
In der nach (19) berechneten Vektorfolge \(\lbrace{\underline{k}_{m} \rbrace}_{m \in \lbrace 0,1,2,...,N-1\rbrace } \) konvergiert für \(N\rightarrow \infty \) und unter der Voraussetzung (15) der Vektor \(\underline{k}_{0}\) gegen einen Grenzwert. Dieser werde wieder mit ''\(\underline{k}\)'' bezeichnet.
Der Grenzwert \(\underline{k}\) ergibt sich auch als Lösung der diskreten algebraischen Matrix - Riccati - Gleichung:
\begin{eqnarray}
\underline{P}=\underline{Q}+\underline{\Phi}^{T}\left[\underline{P} -\underline{P}\ \underline{h}\ \underline{h}^{T}\underline{P}_\cdot \frac{1}{r+\underline{h}^{T}\underline{P}\ \underline{h}}  \right]\underline{\Phi}  
\end{eqnarray}
zu 
\begin{eqnarray}
\underline{k}^{T}=\frac{1}{r+\underline{h}^{T} \underline{P}\ \underline{h}}\cdot \left(\underline{h}^{T}\cdot \underline{P} \cdot \underline{\Phi} \right)
\end{eqnarray}

\end{anm}
Obwohl für die Lösung der diskreten algebraischen Matrix - Riccati - Gleichung schon Lösungsprozeduren bereitgestellt sind (z.B. der Befehl ''DARE'' unter Matlab), empfiehlt es sich im Hinblick auf das betrachtete Problem (iteratives Lösen von linearen Gleichungssystemen), doch besser mit der Matrix - Riccati - Differenzengleichung (20) zu rechnen: alle uns bekannten Lösungsprozeduren für die Lösung der diskreten algebraischen Matrix - Riccati - Gleichung verwenden intern die  Berechnung von Matrixeigenwerten, was bei Gleichungssystemen hoher Ordnung wieder auf Numerikprobleme führen kann. Durch Einsetzen von \(\underline{P}_{m}\) in (21) kann man die Qualität der iterierten Lösung überprüfen. Sei: 
\begin{eqnarray*}
\underline{\Delta}=\underline{Q}+\underline{\Phi}^{T}\left[\underline{P}_{m} -\underline{P}_{m}\ \underline{h}\ \underline{h}^{T}\underline{P}_{m}\cdot \frac{1}{r+\underline{h}^{T}\underline{P}_{m}\ \underline{h}}  \right]\underline{\Phi} -\underline{P}_{m} 
\end{eqnarray*} 
und liege \(||\underline{\Delta}||\) unterhalb einer hinreichend kleinen Schranke, so kann
\(\underline{P}_{m}\) durch \(\underline{P}\) ersetzt werden. Damit berechne man \(\underline{k}^{T}\) entsprechend (22).
\begin{anm}
Mit diesem \(\underline{k}\) - Vektor konvergiert aber auch die Iteration (9):
\begin{eqnarray*}
\underline{\tilde{x}}_{m+1}=\left(\underline{\Phi}-\underline{h}\cdot\underline{k}^{T}\right)\cdot \underline{\tilde{x}}_{m}+\underline{h};\qquad \underline{\tilde{x}}_{m=0}=\underline{\tilde{x}}_{0}
\end{eqnarray*}
gegen einen Fixpunkt.
\end{anm}
Bevor weitere Anmerkungen zu den Matrizen \(\underline{S}\), \(\underline{Q}\) und \(\underline{P}_{m}\) sowie zu ''r'' gemacht werden, sollen die bisherigen Ergebnisse an einem Beispiel demonstriert werden. Wir greifen wieder die \(\underline{A}\) und die \(\underline{\Phi}\) - Matrix sowie den \(\underline{b}\) und \(\underline{h}\) - Vektor von Beispiel (1.3) auf.
\subsection{Weiterführung des Beispiels} 
Mit den recht willkürlichen Werten für \[\underline{c}^{T}=\left(1 \ 4\ 5 \right)\]
und folglich 
\begin{align*}
\underline{Q}=
\begin{pmatrix}
1 & 4 & 5 \\
4 & 16 & 20 \\
5 & 20 & 25
\end{pmatrix}
\end{align*}
sowie mit \[ \underline{S} = \underline{E} \ \text{  -  Einheitsmatrix    und   }r = 0,5 \]
erhält man für N = 100 Iterationen, beginnend mit \(P_{100}\) = Einheitsmatrix mit Hilfe der \newline Gleichung (20)
\begin{align*}
\underline{P}_{99}&\approx
\begin{pmatrix}
1,482.773.652 &  4,200.001.512 &  5,222.175.153 \\
4,200.001.512 & 17,345.911.012 & 19,930.235.464 \\
5,222.175.153 & 19,930.235.464 & 25,132.197.634  
\end{pmatrix}
\end{align*}
\[ \vdots \]
\begin{align*}
\underline{P}_{2}&\approx
\begin{pmatrix}
1,035.281.548 & 4,023.046.495 & 5,014.945.202 \\
4,023.046.495 & 16,060.198.199 & 19,996.629.128 \\
5,014.945.202 & 19,996.629.128 & 25,018.325.346
\end{pmatrix}
\end{align*}

\begin{align*}
\underline{P}_{1}&\approx
\begin{pmatrix}
1,035.281.548 & 4,023.046.495 & 5,014.945.201 \\
4,023.046.495 & 16,060.198.200 & 19,996.629.128 \\
5,014.945.201 & 19,996.629.128 & 25,018.325.346
\end{pmatrix}
\end{align*}

Für die Vektoren \(\underline{k}_{m}^{T}\) ergibt das mit der Gleichung (19):

\begin{align*}
\underline{k}_{99}^{T}&\approx \left(-0,079.698.253\ \ -0,013.703.479\ \ -0,173.741.413\right) \\
\underline{k}_{98}^{T}&\approx \left(-0,007.838.204\quad \ \ \  0,102.143.449\ \ -0,153.895.906\right) \\
\vdots \\
\underline{k}_{1}^{T}&\approx \left(-0,002.607.656\quad \ \ \ 0,106.247.305\ \ -0,151.791.016\right) \\
\underline{k}_{0}^{T}&\approx \left(-0,002.607.656\quad \ \ \ 0,106.247.305\ \ -0,151.791.016\right) \\
\end{align*}
\newline Und kompakter (jetzt mit gerundeten Zahlenwerten):
\newline
\begin{tabular}{|l|r|c|}
\hline
m     & \(P_{m}\qquad \qquad \qquad \qquad \qquad \qquad \)    & \(k_{m}\)    \\
\hline \hline
100 &1,000 \  0,000 \ 0,000\  0,000\ \ 1,000\ \  0,000 \ \  0,000\ \  0,000\ \ 1,000  & --- --- ---  \\

\cline{1-3}
\ 99  &1,483 4,200 5,222 \ 4,200 17,346 19,930 \ 5,222 19,930 25,132 & -0,079.698 \  -0,013.703 \ -0,173.741 \\
\cline{1-3}
\ 98  &1,672 4,629 5,258 \ 4,629 17,168 20,160 \ 5,258 20,160 25,119 & -0,007.838 \ \ 0,102.143 \ -0,153.896 \\
\cline{1-3}
...  & ...\qquad \qquad \qquad \qquad \qquad \qquad \qquad \qquad & ... \qquad \\
\cline{1-3}
\ \ 1  &1,035 4,032 5,015 \ 4,023 16,060 19,997 \ 5,015 19,997 25,018 & -0,002.608 \ \ 0,106.247 \ -0,151.791 \\
\cline{1-3}

\ \ 0  &--- --- --- \qquad \qquad \qquad \qquad \qquad \qquad \qquad & -0,002.608 \ \ 0,106.247 \ -0,151.791 \\
\cline{1-3}
\hline
\end{tabular}

\[\]

Die nachstehende Abbildung 1 zeigt den Verlauf der Komponenten des iterierten \(\underline{k}_{m}\)-Vektors:
\begin{figure}[H]
\begin{center}
\includegraphics[width=16 cm]{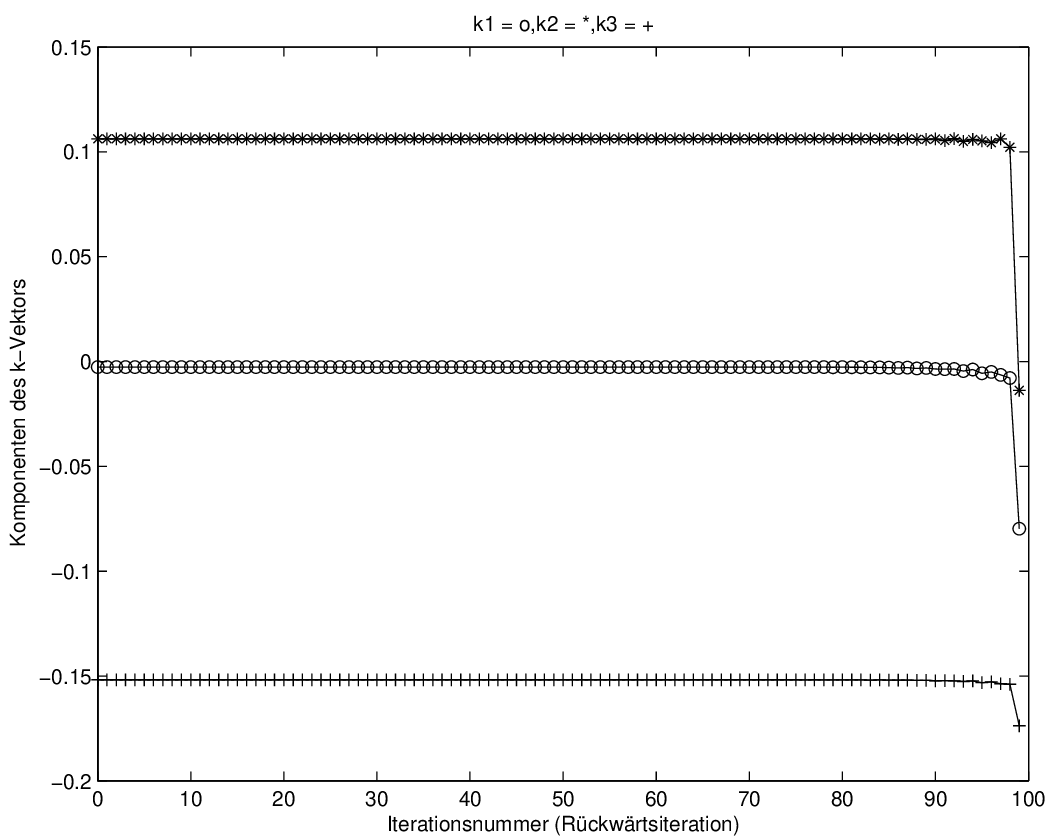}
\caption{Darstellung der Komponenten des \(\underline{k}_{m}\)-Vektors. Die Konvergenz gegen einen Grenzwert ist unter Vorlage der Voraussetzungen von Satz 2.1 und von (15) gesichert.} 
\label{fig:null}
\end{center}
\end{figure}
\(\underline{P}_{1}\) ist eine Näherungslösung von \(\underline{P}\) entsprechend (21):
\begin{eqnarray*}
\underline{P}&=
\begin{pmatrix}
1,035.281.547.600 & 4,023.046.494.700 & 5,014.945.201.200 \\
4,023.046.494.700 & 16,060.198.198.000 & 19,996.629.128.000 \\
5,014.945.201.200 & 19,996.629.128.000 & 25,018.325.346.000
\end{pmatrix}
\end{eqnarray*}
 und \(\underline{k}_{0}\) eine Näherungslösung von (22):
\begin{align*}
\underline{k}^{T}= \left(-0,002.607.655.915\quad \ \ \ 0,106.247.305.105\ \ -0,151.791.015.721\right)
\end{align*}
Mit diesem \(\underline{k}\)-Vektor und dem Startvektor \(\underline{\tilde{x}}_{0}^{T}=(1,1,1)\) liefert die Iteration (9):
\begin{align*}
\underline{\tilde{x}}_{10} = 
\begin{pmatrix}
 \ \  1,140.376.318.225 \\
\ \   1,890.942.195.607\\
\  \  4,622.734.949.616
\end{pmatrix}
\quad
\underline{x}_{10} = 
\begin{pmatrix}
 \ \ 0,758.352.051.190 \\
\ \  1,257.479.544.079 \\
\  \ 3,074.126.036.400
\end{pmatrix}
\end{align*}
\begin{align*}
\underline{\tilde{x}}_{25} = 
\begin{pmatrix}
\ \  1,136.239.118.708 \\
\ \   2,798.412.723.144\\
\  \  3,895.441.810.180
\end{pmatrix}
\quad
\underline{x}_{25} = 
\begin{pmatrix}
 \ \ 0,876.097.560.719 \\
\ \  2,157.717.086.365 \\
\  \ 3,003.581.738.766
\end{pmatrix}
\end{align*}
\begin{align*}
\underline{\tilde{x}}_{100} = 
\begin{pmatrix}
\ \  1,325.309.503.590 \\
\ \   2,650.663.291.458\\
\  \  3,976.119.363.656
\end{pmatrix}
\quad
\underline{x}_{100} = 
\begin{pmatrix}
 \ \ 0,999.954.870.724 \\
\ \  1,999.943.154.231 \\
\  \ 3,000.008.611.948
\end{pmatrix}
\end{align*}
\begin{align*}
\underline{\tilde{x}}_{250} = 
\begin{pmatrix}
\ \  1,325.356.617.751 \\
\ \   2,650.713.325.513\\
\  \  3,976.069.853.246
\end{pmatrix}
\quad
\underline{x}_{250} = 
\begin{pmatrix}
 \ \ 1,000.000.000.002 \\
\ \  2,000.000.000.011 \\
\  \ 2,999.999.999.999
\end{pmatrix}
\end{align*}
(Bei der Berechnung dieser \(\underline{x}\)-Vektoren, wenn statt ''\(\underline{k}\)'' der Vektor ''\(\underline{k}_{0}\)'' verwendet wird, hat man für \(\underline{x}_{10}\) und für \(\underline{x}_{25}\) eine Übereinstimmung in den ersten 9 Nachkommastellen; bei \(\underline{x}_{100}\) und \(\underline{x}_{250}\) jedoch volle Übereinstimmung.)
\newline 
\newline Hierbei beträgt der Restfehler \(\underline{r}_{m}:= \underline{A}\cdot \underline{x}_{m}-\underline{b}\):
\begin{align}
\underline{r}_{100} = 
\begin{pmatrix}
\  \  \ \ 0,901.432 \ 10^{-3} \\
\ \ -3,257.061 \ 10^{-3}\\
\   \ \ \  3,963.873 \ 10^{-3}
\end{pmatrix}
\quad
\underline{r}_{250} = 
\begin{pmatrix}
 \  \ \ -0,619 \ 10^{-9} \\
\quad \ \  0,486 \ 10^{-9}\\
\quad \ \   0,482 \ 10^{-9}
\end{pmatrix}
\end{align}
Für kleine m-Werte sind die Restfehler \(\underline{r}_{m}\) aber noch ''sehr groß''. Die langsame Konvergenz der Folgen  \(\left(\underline{\tilde{x}}_{m}\right)_{m \in \mathbf{N}_{0}}\) und \(\left(\underline{x}_{m}\right)_{m \in \mathbf{N}_{0}}\) wird durch die Eigenwerte von \(\underline{\tilde{\Phi}}\) bedingt. Diese müssen ja nicht mehr vorgegeben werden - die Lösung des Optimierungsproblems garantiert lediglich \(|\tilde{\lambda}_{i}| <1\). Tatsächlich liegen diese mit
\[ |\tilde{\lambda}_{1}|=0,900.828.330.426,\quad |\tilde{\lambda}_{2}|=0,900.828.330.426,\quad |\tilde{\lambda}_{3}|=0,000.430.484.178 \]
zum Teil nahe bei eins.
Die Problematik der gezielten Beeinflussung des Eigenwerttupels (ohne dessen explizite Berechnung) wird in Kapitel 3 aufgegriffen.

\newpage
\subsection{Abschließende Bemerkungen zum Optimierungsverfahren}
Für das in 2.1 vorgestellte Optimierungproblem wird eine Lösung angegeben, die jetzt auch die Berechnung eines linearen Gleichungssystems (1) durch zwei Iterationsverfahren erlaubt:
\begin{enumerate}
\item Iterative Berechnung von \(\underline{P}_{1}\) durch Rückwärtslösen der Matrix-Riccati-Differenzengleichung. Anschließend Berechnung von \(\underline{k}_{0}\).
\item Mit \(\underline{k}=\underline{k}_{0}\) das erweiterte Fixpunktverfahren (9) durchführen. Anschließend Berechnung des Lösungsvektors ''\(\underline{x}\)'' durch Rücktransformation entsprechend (13).
\end{enumerate}
Da jedoch von den Eigenwerten von \(\underline{\tilde{\Phi}}\) nicht mehr gesichert ist, als:
\( \mid \tilde{\lambda}_{i} \mid < 1 \text{ \quad  i} \in \lbrace 1,...,l \leq n \rbrace  \), so kann noch keine Genauigkeitsangabe (a-priori-Fehler und a-posteriori-Fehler) für das erweiterte Iterationsverfahren gemacht werden. Eine solche Angabe gelingt aber, da mit Hilfe des vorgestellten Optimierungsverfahrens auch der betragsgrößte Eigenwert von \(\underline{\tilde{\Phi}}\) abgeschätzt werden kann. Dies wird im nächsten Kapitel ausgeführt.

\section{Ergänzungen zum quadratischen Gütekriterium}

Die Auswahl der Matrizen ''\underline{S}'', des Vektors ''\underline{c}'' und der Zahl ''r'' scheinen zunächst im Hinblick auf die vorgelegte Aufgabenstellung: ''iterative Lösung eines linearen Gleichungssystems'' viele Freiheiten zu lassen. Die Wahl von ''\underline{S}'' erscheint unkritisch, da die Matrix - Riccati - Differenzengleichung unabhängig von ''\underline{S}'' gegen die Lösung der diskreten algebraischen Matrix - Riccati - Gleichung konvergiert. 
\newline Zusätzlich zu der Tatsache, dass für den Spektralradius von \(\underline{\tilde{\Phi}} = \underline{\Phi} - \underline{h}\cdot \underline{k}_{0}^{T}\) gilt: \(\rho\left(\underline{\tilde{\Phi}}\right)<1\),
erlauben die Wahl von ''\(\underline{c}\)'' und ''r'' eine genaue Platzierung der Eigenwerte von 
\( \underline{ \tilde{\Phi} } \) im Innern des Einheitskreises der komplexen Zahlenebene [Hart2]. Da dazu aber auch die Eigenwerte von \( \underline{\Phi} \) zuvor bekannt sein müssen, wird (im Hinblick darauf, dass das hier vorzustellende Iterationsverfahren für Gleichungssysteme hoher Ordnung angewendet werden soll) dieser Aspekt nicht weiterverfolgt.
\newline Immerhin kann der Spektralradius \(\rho\left(\underline{\tilde{\Phi}}\right)\) abgeschätzt werden.

\subsection{Schranken für den Spektralradius \( \rho\left(\underline{\tilde{\Phi}}\right)\). }
Ausgehend von Iteration (14) mit \(\underline{\Phi} \rightarrow \frac{1}{w}\cdot \underline{\Phi}=:\underline{\Phi}_{w}\) und \(\underline{h}\rightarrow \frac{1}{w} \cdot \underline{h}=:\underline{h}_{w}\)\qquad    (\(w>0\))
\newline und
\begin{eqnarray}
\underline{x}_{m+1}= \underline{\Phi}_{w}\ \cdot \underline{x}_{m}+\underline{h}_{w}\cdot u_{m}\ ;\qquad \underline{x}_{m=0}=\underline{x}_{0} 
\end{eqnarray}
werde eine Folge \( \lbrace u_{m} \rbrace_{m \in \mathbf{N}_{0}} \) gesucht, sodass mit 
\( u_{m}=-\underline{k}_{m}^{T} \cdot \underline{x}_{m} \)  das quadratische Gütekriterium (18) minimal wird. Sei \(\underline{k}_{w}^{*}\) der ''optimale Vektor'' der sich als 
\( \lim_{N\rightarrow \infty} \underline{k}_{0} \) oder als Lösung der diskreten algebraischen Matrix - Riccati - Gleichung ergibt.
Mit diesem Vektor \(\underline{k}_{w}^{*}\) gilt 
\begin{eqnarray}
\rho \left(\underline{\Phi}_{w} + \underline{h}_{w}\cdot \underline{k}_{w}^{*T}\right) < 1
\end{eqnarray}
Bezeichne nunmehr \(\lambda_{i}^{*}\) einen Eigenwert von \( \frac{1}{w}\cdot \underline{\Phi} + \frac{1}{w}\cdot \underline{h}\cdot \underline{k}_{w}^{*T} \) 
\newline und \(\tilde{\lambda}_{i}\) (wie oben) einen Eigenwert von \(\underline{\tilde{\Phi}_{w}}:= \underline{\Phi} - \underline{h}\cdot \underline{k}_{w}^{*T} \) sowie \(\underline{E}\) die Einheitsmatrix, so folgt aus
\begin{eqnarray*}
0 = det\left( \lambda_{i}^{*} \cdot \underline{E} -  \left( \frac{1}{w}\cdot \underline{\Phi} + \frac{1}{w}\cdot \underline{h}\cdot \underline{k}_{w}^{*T}\right) \right)=det\left( w\cdot\lambda_{i}^{*} \cdot \underline{E} -  \left( \underline{\Phi} + \underline{h}\cdot \underline{k}_{w}^{*T}\right) \right) 
\end{eqnarray*}
endlich \(  \tilde{\lambda}_{i}= w\cdot\lambda_{i}^{*} \). Mit diesem \( \underline{k}_{w}^{*} \) gilt: 
\begin{eqnarray}
\rho\left(\underline{\tilde{\Phi}}\right)<w
\end{eqnarray}  
Das so erweiterte Iterationsverfahren wird mit den Rechnungen des nachfolgenden Unterpunkts vorgestellt. Es wird wieder das unter 1.3 vorgestellte Rechenbeispiel aufgegriffen; diesmal aber mit 9 bzw. 12 Nachkommastellen.

\subsection{Rechenbeispiele zum erweiterten Iterationsverfahren}
Mit
\begin{align*}
\underline{A} = 
\begin{pmatrix}
59 & -63 & -2 \\
29 & 42 & 51 \\
36 & 31 & -67
\end{pmatrix}
\quad \text{  und } \qquad \underline{b}=
\begin{pmatrix}
-73 \\
266\\
-103
\end{pmatrix}
\end{align*}
berechnen sich die \(\underline{\Phi}\) - Matrix und der \(\underline{h}\) - Vektor wie in 1.3.
\newline Mit \(w = \frac{1}{10}\) und \( \underline{\Phi}_{w}= \frac{1}{w}\cdot \underline{\Phi} \),  
\( \underline{h}_{w}= \frac{1}{w}\cdot \underline{h} \) lauten
\begin{align*}
\underline{\Phi}_{w}= 
\begin{pmatrix}
\ \ 0,000.000.000 & 10,677.966.102 & \quad \  0,338.983.051 \\
-6,904.761.905 & \ \ 0,000.000.000 &-12,142.857.143 \\
\ \ 5,373.134.328 & 4,626.865.672 & 0,000.000.000
\end{pmatrix} 
\text{  und  }
\quad \underline{h}_{w}=
\begin{pmatrix}
&-12,372.881.355.932 \\
&\ 63,333.333.333.333\\
&\ 15,373.134.328.358
\end{pmatrix}
\end{align*}
Unter Verwendung von \(\underline{S}\), \(\underline{c}\) und ''r'' wie in 2.2 erhält man für N = 10 und der Einheitsmatrix für \(\underline{P}_{10}\) mit Gleichung (20):
\begin{align*}
\underline{P}_{9}&\approx
\begin{pmatrix}
48,959.414.201 &  29,945.482.078 & 26,524.385.209 \\
23,945.482.078 & 150,581.701.310 & 12,904.368.225 \\
26,524.385.209 & 12,904.368.225 & 36,708.746.563  
\end{pmatrix}
\end{align*}
\(\qquad \qquad \qquad \qquad \qquad \qquad \qquad \qquad\qquad \qquad \qquad  \vdots\)
\begin{align*}
\underline{P}_{1}&\approx
\begin{pmatrix}
6.195,608.211.100 & 30.591,456.824.000 & -774,625.645.620 \\
30.591,456.824.000 & 181.987,753.000.000 & -7.760,651.406.800 \\
-774,625.645.620 & -7.760,651.406.800 & 622,611.391.310
\end{pmatrix}
\end{align*}
Für die Vektoren \(\underline{k}^{T}_{m}\) ergibt das mit Gleichung (19):
\begin{align*}
\underline{k}_{9}^{T}&\approx \left(-0,080.594.648.281 \quad -0,013.857.607.126 \quad -0,175.695.545.308           \right) \\
\vdots \\
\underline{k}_{0}^{T}&\approx \left(-0,117.736.653.973\quad \ \ \ 0,026.270.720.999\ \ -0,199.437.647.693\right) 
\end{align*}
Mit \(\underline{k}^{*}_{w}=\underline{k}_{0}\) folgt sodann mit:
\begin{align*}
\underline{\tilde{\Phi}}_{w}:= \underline{\Phi} - \underline{h}\cdot \underline{k}_{w}^{*T} 
\end{align*}
für  
\begin{align*}
\tilde{\underline{x}}_{m+1}= \underline{\tilde{\Phi}}_{w} \cdot \tilde{\underline{x}}_{m}+\underline{h} 
\end{align*}
und für
\begin{align*}
\underline{x}_{m}= \frac{1}{1-\underline{k}_{w}^{*T}\cdot \underline{\tilde{x}}_{m}}\cdot \underline{\tilde{x}}_{m} \qquad \text{ sowie } \qquad \underline{r}_{m}=\underline{A}\cdot \underline{x}_{m}-\underline{b}
\end{align*}
und den Startvektor  \(\ \underline{\tilde{x}}_{0}^{T} = (1, 1, 1) \) die Vektoren:
\begin{align*}
m+1 = 5: \quad \quad
\underline{\tilde{x}}_{5} = 
\begin{pmatrix}
 \ \  2,971.843.010.215 \\
\ \   5,943.680.457.527\\
\  \  8,915.516.643.666
\end{pmatrix}
\quad 
\underline{x}_{5} = 
\begin{pmatrix}
 \ \ 1,000.001.061.532 \\
\ \  2,000.000.251.191 \\
\  \  2,999.999.016.477
\end{pmatrix}
\quad
\underline{r}_{5} = 
\begin{pmatrix}
 \  \  \  \ 0,000.048.722 \\
\ \ -0,000.008.825   \\
\  \ \ \  0,000.111.898
\end{pmatrix}
\end{align*}

\begin{align*}
m+1 = 10: \quad
\underline{\tilde{x}}_{10} = 
\begin{pmatrix}
 \ \  2,971.840.542.360 \\
\ \   5,943.681.084.719\\
\  \  8,915.521.627.079
\end{pmatrix}
\quad
\underline{x}_{10} = 
\begin{pmatrix}
 \ \ 1,000.000.000.000 \\
\ \  2,000.000.000.000 \\
\  \ 3,000.000.000.000
\end{pmatrix}
\quad
\underline{r}_{10} = 
\begin{pmatrix}
 -1,421.085 \ 10^{-14} \\
\ \ 0,000.000 \  10^{-14}   \\
-2,842.171 \ 10^{-14}
\end{pmatrix}
\end{align*}
Die Mächtigkeit des erweiterten Iterationsverfahrens zeigt sich erst bei Erniedrigung des Parameters ''w'' !
\[\] Seien \(w = \frac{1}{100}\) und die anderen Parameter unverändert, dann erhalten wir nach N = 5 Iterationsschritten (beginnend mit \(\underline{P}_{5}\) = Einheitsmatrix):
\begin{align*}
\underline{P}_{1}&\approx
\begin{pmatrix}
\ \ 4,877.396.786 \ 10^{7} & 2,975.175.473 \  10^{8} & -1,334.071.540 \ 10^{7} \\
\ \ 2,975.175.473 \  10^{8} & 1,818.776.612 \ 10^{9} & -8.187.822.688 \ 10^{7} \\
-1,334.071.540 \ 10^{7}&-8.187.822.688 \ 10^{7} & \  3,712.634.543 \ 10^{6}
\end{pmatrix}
\end{align*}

\begin{align*}
\underline{k}_{0}^{T}&\approx \left(-0,117.898.057.028\quad \ \ \ 0,025.378.263.771\ \ -0,199.405.203.705\right) 
\end{align*}
und nach M = m+1 = 5 Schritten für die Fixpunktiteration
\begin{align*}
M = 5: \quad
\underline{\tilde{x}}_{5} = 
\begin{pmatrix}
 \ \  2,988.260.385.409 \\
\ \   5,976.520.770.812\\
\  \  8,964.781.156.214
\end{pmatrix}
\quad
\underline{x}_{5} = 
\begin{pmatrix}
 \ \ 1,000.000.000.000 \\
\ \  2,000.000.000.000 \\
\  \ 3,000.000.000.000
\end{pmatrix}
\quad
\underline{r}_{5} = 
\begin{pmatrix}
\ \  4,908 \  10^{-11} \\
-8,868 \ 10^{-12}   \\
\ \ 1,112 \ 10^{-10}
\end{pmatrix}
\end{align*}
\[ \]
Seien jetzt noch \(w = \frac{1}{10.000}\ \),  N = 2 und M = m+1 = 4, dann erhalten wir (allen anderen Parameter unverändert):
\begin{align*}
\underline{P}_{1}&\approx
\begin{pmatrix}
\ \ 2,061.646.012 \  10^{14} & \ 1,260.667.803 \ 10^{15} & -5,678.165.532 \ 10^{13} \\
\ \ 1,260.667.803 \ 10^{15} & \ 7,708.847.447 \ 10^{15} & -3,472.169.309  \  10^{14} \\
-5,678.165.532 \  10^{13} & -3,472.169.309 \ 10^{14} & \ 1,563.938.679   \ 10^{13}
\end{pmatrix}
\end{align*}
\begin{align*}
\underline{k}_{0}^{T}&\approx \left(-0,117.899.526.977\quad \ \ \ 0,025.369.284.637\ \ -0,199.404.800.048\right) 
\end{align*}

\begin{align*}
M = 4: \quad
\underline{\tilde{x}}_{4} = 
\begin{pmatrix}
 \ \  2,988.423.068.799 \\
\ \   5,976.846.137.597\\
\  \  8,965.269.206.396
\end{pmatrix}
\quad
\underline{x}_{4} = 
\begin{pmatrix}
 \ \ 1,000.000.000.000 \\
\ \  2,000.000.000.000 \\
\  \ 3,000.000.000.000
\end{pmatrix}
\quad
\underline{r}_{4} = 
\begin{pmatrix}
\ \  7,958  \ 10^{-13} \\
-1,705 \ 10^{-13}   \\
\ \ 1,691 \ 10^{-12}
\end{pmatrix}
\end{align*}

\subsection {Algorithmische Formulierung des erweiterten Iterationsverfahrens}
Vorbereitend für die algorithmische Formulierung werden die folgenden Definitionen eingeführt:
\begin{defi}{Die Funktion Ricc}
\begin{align*}
Ricc:\qquad \mathbf{R}^{n \times n} \longrightarrow \mathbf{R}^{n \times n} \qquad \qquad \qquad \qquad \qquad \qquad \qquad \qquad \qquad \qquad \qquad \qquad \qquad \\
\underline{P} \longmapsto Ricc\left(\underline{P}\right)= \underline{Q}+\underline{\Phi}^{T}_{w}\cdot \left( \underline{P}-\underline{P} \cdot \left(\underline{h}_{w}\cdot\underline{h}^{T}_{w}\right) \cdot \underline{P}\cdot 
\frac{1}{r + \underline{h}_{w}\cdot \underline{P} \cdot\underline{h}^{T}_{w}}                        \right)
\cdot \underline{\Phi}_{w} \quad  
\end{align*}
\end{defi}
\begin{defi}{Die Funktion \(\underline{F}\)}
\begin{align*}
F:\qquad \mathbf{R}^{n} \longrightarrow \mathbf{R}^{n} \qquad \qquad \qquad \qquad \qquad \qquad \qquad \qquad \qquad \qquad \qquad \qquad \qquad \\
\underline{x} \longmapsto \underline{F}(\underline{x})=\left( \underline{\Phi}_{w}-\underline{h}_{w}\cdot \underline{k}_{0}^{T}\right)\cdot\underline{x} \qquad \qquad \qquad \qquad \qquad \qquad \qquad \qquad  \ \end{align*}
\end{defi}
Nun zum Algorithmus.
\newpage
\begin{enumerate}
\item Zum linearen Gleichungssystem \(\underline{A}\cdot \underline{x}=\underline{b}\ \):   
\(\ \underline{\Phi}\) und \(\underline{h}\) berechnen.
\item \(0 < w <1\ \), N und M wählen.
\item \(\underline{\Phi}_{w}:=\frac{1}{w} \cdot \underline{\Phi}\) und \(\underline{h}_{w}:=\frac{1}{w}\cdot\underline{h}\ \)  berechnen.
\item Mit Startwert \(\underline{P}_{N}=\underline{E} \ \)- Einheitsmatrix in N - Iterationsschritten mittels Funktion ''Ricc'' durch
\begin{align*}
\underline{P}_{n-1}=Ricc \left( \underline{P}_{n} \right)
\end{align*}
den Näherungswert \(\underline{P}_{1}\) berechnen.
\item Den Vektor \(\underline{k}_{0}\) berechnen:
\begin{align*}
\underline{k}_{0}^{T}=\frac{1}{r+\underline{h}_{w}^{T}\cdot \underline{P}_{1}\cdot \underline{h}_{w}}\cdot \left(\underline{h}_{w}^{T}\cdot \underline{P}_{1} \cdot \underline{\Phi}_{w} \right) 
\end{align*}
\item Mit beliebigem Startwert \(\underline{\tilde{x}}_{0}\) in M - Iterationsschritten (\(M > 1\)) mittels der Funktion ''\(\underline{F}\) '' durch
\begin{align*}
 \underline{\tilde{x}}_{m+1}= \underline{F} \left( \underline{\tilde{x}}_{m}\right) +\underline{h}
\end{align*}
den Näherungswert \(\underline{\tilde{x}}_{M}\) berechnen.
\item Den Näherungswert \(\underline{x}_{M}\) gemäß Vorschrift:
\begin{align*}
\underline{x}_{M}= \left(\frac{1}{1-\underline{k}_{0}^{T} \cdot \underline{\tilde{x}}}_{M}\right) \cdot \underline{\tilde{x}}_{M}
\end{align*}
und 
\begin{align*}
\underline{r}_{M}=\underline{A}\cdot \underline{x}_{M}- \underline{b}
\end{align*}
berechnen.
\end{enumerate}


\printindex
\end{document}